\documentclass[11pt,reqno]{amsart}
\usepackage{stmaryrd}
\usepackage{xcolor}

\usepackage{epsfig}
\usepackage{enumerate}
\usepackage{amsmath}
\usepackage{amssymb}
\usepackage{amscd}
\usepackage{graphicx}
\usepackage{lineno}
\usepackage{pstricks}
\usepackage{tocvsec2}
\usepackage{mathtools}
\usepackage{bm}
\allowdisplaybreaks[4]

\usepackage[
backref,
colorlinks,
citecolor=red,
anchorcolor=black,
linkcolor=blue
]{hyperref}

\usepackage{amsfonts}
\usepackage[
top=35mm,
bottom=35mm,
left=30mm,
right=30mm
]{geometry}

\usepackage{mathrsfs}

\makeatletter
\@namedef{subjclassname@2020}{%
	\textup{2020} Mathematics Subject Classification}
\makeatother

\newcommand{\la}{\Lambda^{(s)}_{k}}
\newcommand{\J}{J_{k,s}(\bm{\lambda})}

\newcommand{\lam}{\bm{\lambda}}
\newcommand{\x}{\bm{x}}

\makeatletter
\NewCommandCopy\@@pmod\pmod
\DeclareRobustCommand{\pmod}{\@ifstar\@pmods\@@pmod}
\def\@pmods#1{\mkern4mu({\operator@font mod}\mkern 6mu#1)}
\makeatother

\newtheorem{lemma}{Lemma}[section]

\newtheorem{theorem}{Theorem}[section]

\renewcommand*{\backref}[1]{}

\numberwithin{equation}{section}

\begin{document}
	
	\title[ Random Diophantine inequalities in few variables]{ Random Diophantine inequalities in few variables}

\author{Yu-Chen Sun}
\address[Yu-Chen Sun]
{School of Mathematics, University of Bristol, Bristol, BS8 1UG, England}
\email{yuchensun93@163.com}
\author{Xuanzhe Wang}
\address[Xuanzhe Wang]
{School of Mathematical Sciences, University of Science and Technology of China, Hefei, 230026, China}
\email{wxzdsgdddsg@mail.ustc.edu.cn}
\author{Lilu Zhao}
\address[Lilu Zhao]
{School of Mathematical Sciences, University of Science and Technology of China, Hefei 230026, China}
\email{zhaolilu@ustc.edu.cn}

\begin{abstract}We prove that for almost all additive Diophantine inequalities of degree $k$ in at least $k+1$ variables the expected asymptotic formula for the density of solutions holds true. This improves the result of Br\"udern and Dietmann \cite{BD}.
\end{abstract}

\maketitle

	\section{Introduction}
	
	The study of additive Diophantine inequalities is one of the important topics in Diophantine approximations. For the classical work on the Diophantine approximation, one may refer to  Hardy-Littlewood \cite{HL}, 
	Davenport and Heilbronn \cite{DH} and Schmidt \cite{Sch}. For Diophantine inequalities involving prime variables, one may refer to Vaughan \cite{Va1974}. Let $k$ and $s$ be natural numbers with $s>k\geqslant 2$. Let
	\begin{align*}
		F_{\lam}(\mathbf{x})=\lambda_1x_1^k+\lambda_2x_2^k+\cdots+\lambda_sx_s^k,
	\end{align*}
where $\lambda_1,\ldots,\lambda_s$ are nonzero real numbers. A natural assumption is that at least one ratio $\lambda_i/\lambda_j$ is irrational. A further assumption is that the coefficients must have mixed signs when $k$ is even. Given $\tau>0$, an important problem is to seek nonzero integer vectors $\mathbf{x}\in \mathbb{Z}^{s}$ such that
$$|F_{\lam}(\mathbf{x})|<\tau.$$
We introduce the counting function \begin{align*}
		N_{\bm{\lambda}}^{(k)}(P;\tau)=\#\left\lbrace \x\in \mathbb{Z}^s:\left|
		\x\right|_{\infty}\leqslant P,\;|F_{\lam}(\x)|<\tau\right\rbrace,
	\end{align*}
where $|\mathbf{x}|_{\infty}=\max_{1\leqslant j\leqslant s} |x_j|$. The expected asymptotic formula is 
\begin{align}\label{asymptotic}N_{\bm{\lambda}}^{(k)}(P;\tau)\sim 2\tau J_{k,s}(\bm{\lambda}) P^{s-k},\end{align}
where the singular integral $J_{k,s}(\bm{\lambda})$ describes the real density of the zero set in the unit box. 
Motivated by the work of Bentkus and G\"otze \cite{BG}, Freeman \cite{F1,F2} established the asymptotic formula \eqref{asymptotic} 
when $s\geqslant 2^k+1$. The number of variables involved in Freeman's asymptotic result can be considerably improved since the work of \cite{W2,W3}. One can refer to \cite{BDG,W4} for the development after \cite{W3}. 
When $k=2$, the celebrated work
	of Margulis \cite{Ma} answers a long-standing conjecture of Oppenheim \cite{Oppenheim}. Wooley \cite{W2003} obtained a lower bound
 $$N_{\bm{\lambda}}^{(k)}(P;\tau)\gg \tau J_{k,s}(\bm{\lambda}) P^{s-k}$$
 when $s\geqslant k(\log k+\log\log k+2+o(1))$.

We write $\Lambda^{(s)}=(\mathbb{R}\setminus\{0\})^s$. Then we define
 $\Lambda_k^{(s)}=\Lambda^{(s)}$ when $k$ is odd, and
	\begin{align*}
		\Lambda_k^{(s)}
		=
		\left\{
		\boldsymbol{\lambda}\in\Lambda^{(s)}:
		\lambda_i\lambda_j<0
		\text{ for some }1\leqslant i<j\leqslant s
		\right\}
	\end{align*}
	when $k$ is even. Br\"udern and Dietmann \cite{BD} proved that \eqref{asymptotic} holds for almost all $\lam$ in the sense of Lebesgue measure whenever $s>2k$. Precisely, they proved the following.

	\begin{theorem}[Br\"udern and Dietmann]\label{thm1}Let $s,k$ be integers with $s>2k$ and $k\geqslant 2$. Suppose that $\delta=8^{-2k}$.
		Then for almost all $\bm{\lambda}\in\la$,
		there exists a number
		$P_0=P_0(\bm{\lambda},k)$
		and a positive real number
		$J_{k,s}(\bm{\lambda})$
		such that the inequality
		\begin{align*}
			\left|
			N_{\bm{\lambda}}^{(k)}(P;\tau)
			-2\tau J_{k,s}(\bm{\lambda}) P^{s-k}
			\right|
			<P^{s-k-\delta}
		\end{align*}
		holds for all $P\geqslant P_0$
		and $0<\tau\leqslant 1$.
	\end{theorem}
	
Another interesting direction on the exceptional sets for Diophantine inequalities was developed by Parsell and Wooley \cite{PW}, who studied the exceptional real numbers $\xi$ such that 
$$\left|\lambda_1x_1^k+\lambda_2x_2^k+\cdots+\lambda_sx_s^k-\xi\right|<\tau$$
has no integer solutions.

The aim of this paper is to prove the following result which improves Theorem \ref{thm1}.

	\begin{theorem}\label{thm2}Let $s,k$ be integers with $s>k\geqslant 2$. Suppose that $\delta= \frac{1}{20k}$.
		Then for almost all $\bm{\lambda}\in\la$,
		there exists a number
		$P_0=P_0(\bm{\lambda},k)$
		and a positive real number
		$J_{k,s}(\bm{\lambda})$
		such that the inequality
		\begin{align*}
			\left|
			N_{\bm{\lambda}}^{(k)}(P;\tau)
			-2\tau J_{k,s}(\bm{\lambda}) P^{s-k}
			\right|
			< P^{s-k-\delta}
		\end{align*}
		holds for all $P\geqslant P_0$
		and $0<\tau\leqslant 1$.
	\end{theorem}
	
Let us briefly discuss the idea of the proof. We follow the general framework of Br\"udern and Dietmann \cite{BD}, using a triangular weight, the Davenport--Heilbronn method, and a Borel--Cantelli argument followed by removal of the weight. The main difference lies in the treatment of the minor arc contribution.

Let $I(\lam;\mathfrak{m})$, $\la(C)$ and $f(\alpha)$ be defined as in \eqref{defineID}, \eqref{definelC} and \eqref{definef}, respectively, and let $K_{\eta}$ be the kernel introduced in Section~\ref{sec:notation}.
To estimate the minor arc contribution, Br\"udern and Dietmann bound the first moment by an integral of the absolute value of the product of Weyl sums,
\[
\int_{\la(C)}|I(\lam;\mathfrak{m})|\,d\lam
\leqslant \int_{\la(C)}\int_{\mathfrak{m}}
\prod_{j=1}^{s}|f(\lambda_j\alpha)|K_{\eta}(\alpha)\,d\alpha\,d\lam.
\]
They combine a mean value estimate involving $2k$ factors with a saving from an additional factor, and use Weyl's inequality to bound the measure of the set on which all the Weyl sums are large; see \cite[Lemmas 4 and 6]{BD}.

We instead estimate the second moment of the minor arc integral by 
\begin{align*}
\int_{\la(C)}|I(\lam;\mathfrak{m})|^2\,d\lam & =\int_{\la(C)}\int_{\mathfrak{m}}\int_{\mathfrak{m}}\prod_{j=1}^{s} f(\lambda_j\alpha)\overline{f(\lambda_j\beta)}\,d\alpha\,d\beta\,d\lam \\
&=\int_{\mathfrak{m}}\int_{\mathfrak{m}} \left(\sum_{|x|\leq P}\sum_{|y|\leq P}\int_{-C}^{C}e(\lambda(\alpha x^k - \beta y^k))\,d\lambda \right)^s\,d\alpha\,d\beta\
\end{align*}
Note that the integral from $-C$ to $C$ is a Dirichlet kernel, which is bounded by 
\begin{align*}
\sum_{|x|,|y|\leqslant P}
\min\left(1,\frac{1}{|\alpha x^k-\beta y^k|}\right).
\end{align*}
Thus, the problem is reduced to estimating the contribution of pairs for which $\alpha x^k$ and $\beta y^k$ are close.

The elementary spacing estimate in Lemma \ref{lem1}, inspired by the large sieve inequality, then yields the second moment bound in Lemma \ref{minor arc}, without using pointwise minor arc estimates for the Weyl sums.

The paper is organised as follows. In Section~\ref{sec:notation}, we introduce the notation and state the major and minor arc estimates. Section~\ref{sec:major} establishes the major arc estimate. In Section~\ref{sec:minor}, we prove the spacing lemma and use it to bound the second moment of the minor arc contribution. In Section~\ref{sec:proof}, we apply the Borel--Cantelli lemma and remove the weights to complete the proof of Theorem~\ref{thm2}. 

	\section{Notation}\label{sec:notation}
	In this section, we shall introduce some notation and provide the asymptotic result on the major arcs. As usual, we write $e(\alpha)=e^{2\pi i\alpha}$.  Throughout this paper, in Vinogradov's notation, the implicit constants in $\ll$ and $\gg$ may depend on $k,s,C$, but are independent of $\eta$ and $\tau$.
	Let us define the Fej\'er kernel as follows.
	\begin{align*}
		K_{\eta}(\alpha)=\eta \left(\frac{\sin(\pi \eta \alpha)}{\pi \eta \alpha}\right)^2,\;\;\;W_{\eta}(t)=\frac{1}{\eta}(\mathbf{1}_{(0,\eta)} * \mathbf{1}_{(- \eta,0)})(t)=\max\left(0,\ 1-\frac{\left|t\right|}{\eta}\right)	,
	\end{align*} 
	where $0<\eta\leqslant 2$. One has
	\begin{align}\label{propK1}
		K_{\eta}(\alpha)\ll \eta \min\left(1,\frac{1}{(\eta \alpha)^2}\right)\ll \frac{\eta}{1+(\eta\alpha)^2}.
	\end{align}
In particular, 
	\begin{align}\label{propK0}
		K_{\eta}(\alpha)\ll \eta.
	\end{align}
    The Fourier transform gives
    \begin{align}\label{FT}
    	\int_{\mathbb{R}}K_{\eta}(\alpha)e(\alpha t)d\alpha=W_{\eta}(t).
    \end{align}
    In particular, we have
    \begin{align}\label{PropK2}
    	\int_{\mathbb{R}}K_{\eta}(\alpha)d\alpha=1.
    \end{align}      
	We define
	\begin{align*}
		I(P,\eta,\lam)=\sum_{\left|\x\right|_{\infty}\leqslant P}W_{\eta}(F_{\lam}(\x)).
	\end{align*}
	Let
	\begin{align}\label{definef}
		f(\alpha)=\sum_{\left|x\right|\leqslant P}e(\alpha x^k).
	\end{align}
	By \eqref{FT}, we have
	\begin{align*}
		I(P,\eta,\lam)=\int_{\mathbb{R}}f(\lambda_1\alpha)f(\lambda_2\alpha)\cdots f(\lambda_s\alpha)K_{\eta}(\alpha)d\alpha.
	\end{align*}

	Let $C\geqslant 1$ be a fixed real number. We define
	\begin{align}\label{definelC}
		\la(C)=\left\lbrace \lam\in \la:1/C\leqslant \left|\lambda_i\right|\leqslant C(1\leqslant i\leqslant s)\right\rbrace 
	\end{align}
	 and write  
\begin{align}\label{defineQ}Q=(4kC)^{-1}P^{1-k}.\end{align}
We now define the major arc
		$$\mathfrak{M}=[-Q,Q]$$ 
    and the minor arc 
    $$\mathfrak{m}=\mathbb{R}\setminus\mathfrak{M}.$$
    	For $\mathfrak{D}\in \{\mathfrak{M},\,\mathfrak{m}\}$, we write
	\begin{align}\label{defineID}
		I(\lam;\mathfrak{D})=\int_{\mathfrak{D}}f(\lambda_1\alpha)f(\lambda_2\alpha)\cdots f(\lambda_s\alpha)K_{\eta}(\alpha)d\alpha.
	\end{align}
    Then we have	\begin{align*}
		I(\lam;\mathfrak{D})=I(\lam;\mathfrak{M})+I(\lam;\mathfrak{m}).
	\end{align*}

We present the major arc estimate as follows. 
        \begin{lemma}\label{major arc}
    	Suppose $0<\eta\leqslant 2$. Let $C\geqslant 1$ and let $\lam\in \la(C)$. There exists a continuous function $J_{k,s}:\, \Lambda_{k}^{(s)}\rightarrow \mathbb{R}^{+}$ such that
    	\begin{align*}
    		I(\lam;\mathfrak{M})=\eta J_{k,s}(\bm{\lambda}) P^{s-k}+O(\eta P^{s-k-\frac{1}{k}}).
    	\end{align*}
    \end{lemma}

The next is the minor arc estimate.
    \begin{lemma}\label{minor arc}
    	Suppose that $0< \eta\leqslant 2$. For any fixed $C>1$, we have
    	\begin{align*}
    		\int_{\la(C)}\left| I\left(\lam ;\mathfrak{m}\right)\right|^2d\lam \ll  1+\eta P^{s-k}(\log P)^{4s}+\eta^2P^{2s-2k-\frac{2s-2k}{k}}(\log P)^{4s}.
    	\end{align*}
    \end{lemma}  
    
    Before proceeding to the next section, we give a remark. Define 
    	\begin{align*}
		F(\alpha)=\sum_{\left|x\right|\leqslant P}\xi_{x}e(\alpha x^k)
	\end{align*}where the sequence $\xi_{x}$ satisfies $|\xi_{x}|\leqslant 1$ for all $x$. 
Let
\begin{align*}\mathcal{I}_{\xi}(L)=\int_{\la(C)}\left|\int_{LP^{-k}}^{+\infty}F(\lambda_1\alpha)F(\lambda_2\alpha)\cdots F(\lambda_s\alpha)K_{\eta}(\alpha)d\alpha)\right|^2d\lam.
	\end{align*}
In fact, the proof of Lemma \ref{minor arc} implies the following
    	\begin{align}\label{forprime}
    		\mathcal{I}_{\xi}(L) \ll  1+\eta P^{s-k}(\log P)^{4s}+\eta^2P^{2s-2k}L^{-\frac{2s-2k}{k}}(\log P)^{4s}.
    	\end{align}
    One can use \eqref{forprime} to establish a variant of Theorem \ref{thm2} in $k+1$ prime variables.
    
   \section{The major arcs estimate}\label{sec:major}
   We give the proof of Lemma \ref{major arc} in this section.
\begin{proof}[Proof of Lemma \ref{major arc}] Since the proof of the major arc estimate is standard, we present it briefly.
	Let
	\begin{align*}
		v(\alpha)=\int_{-P}^{P}e(\alpha t^k)dt
	\end{align*}
and 	\begin{align*}
		V_{\lam}(\alpha)=\prod_{i=1}^{s}v(\lambda_i\alpha).
	\end{align*}
	The partial integration yields
	\begin{align}\label{esv}
		v(\alpha)\ll P(1+P^k\left|\alpha\right|)^{-\frac{1}{k}}.
	\end{align}
Recalling \eqref{defineQ} and applying \cite[Theorem 4.1]{Va}, we have
	\begin{align}\label{ap}
		f(\lambda_i\alpha)=v(\lambda_i\alpha)+O(1),\;1\leqslant i\leqslant s,
	\end{align}
    when $\lam\in \la(C)$ and $\alpha\in \mathfrak{M}$.  Then \eqref{esv} and \eqref{ap} give
    \begin{align*}
    	f(\lambda_1\alpha)\cdots f(\lambda_s\alpha)-V_{\lam}(\alpha)\ll P^{s-1}\left(1+P^k\left|\alpha\right|\right)^{-\frac{s-1}{k}},
    \end{align*}
    and thus,
    \begin{align*}
    	I(\lam;\mathfrak{M})-\int_{\mathfrak{M}}V_{\lam}(\alpha)K_{\eta}(\alpha)d\alpha\ll \int_{\mathfrak{M}}P^{s-1}\left(1+P^k\left|\alpha\right|\right)^{-\frac{s-1}{k}}K_{\eta}(\alpha)d\alpha .
    \end{align*}
    Using \eqref{propK0}, we obtain
        \begin{align}\label{err1}
    	I(\lam;\mathfrak{M})-\int_{\mathfrak{M}}V_{\lam}(\alpha)K_{\eta}(\alpha)d\alpha\ll  \eta P^{s-1-k}\log P,
    \end{align}
   where the factor $\log P$ arises (only) in the case $s=k+1$.

    Using Taylor's expansion, we have
    \begin{align*}
    	\left|K_{\eta}(\alpha)-\eta\right|\ll \eta^3\alpha^2.
    \end{align*}
Then we deduce that
    \begin{align}\label{err2}
    	\int_{\mathfrak{M}}V_{\lam}(\alpha)K_{\eta}(\alpha)d\alpha -\eta\int_{\mathfrak{M}}V_{\lam}(\alpha)d\alpha\ll \int_{\mathfrak{M}}\eta^3|\alpha|^2P^{s}\left(1+P^k\left|\alpha\right|\right)^{-\frac{s}{k}}d\alpha\ll \eta^3P^{s-3k+2}.
    \end{align}
    Note that for $k\geqslant 2$ and $\eta\leqslant 2$, one has
    $$\eta^3P^{s-3k+2}\ll \eta P^{s-k-1}.$$ 
By \eqref{esv}, we have
    \begin{align}\label{err3}
    	\eta\int_{\mathfrak{M}} V_{\lam}(\alpha)d\alpha-\eta\int_{\mathbb{R}} V_{\lam}(\alpha)d\alpha\ll \eta P^{s-k-\frac{s-k}{k}}.
    \end{align}
    Combining \eqref{err1}, \eqref{err2} and \eqref{err3},  we obtain
            \begin{align}\label{err4}
    	I(\lam;\mathfrak{M})-\eta\int_{\mathbb{R}}V_{\lam}(\alpha)d\alpha \ll  \eta P^{s-1-k}\log P+\eta P^{s-k-\frac{s-k}{k}}.
    \end{align}

    We define
    \begin{align*}
    	J_{k,s}(\bm{\lambda})=\int_{\mathbb{R}}\prod_{i=1}^{s}\left(\int_{-1}^{1}e(\alpha\lambda_i t^k)dt\right)d\alpha.
    \end{align*}
    The property of $J_{k,s}(\bm{\lambda})$ can be found in \cite[section 5]{BD}. In particular, it is positive. By a series of changing variables, we can obtain
            \begin{align}\label{err5}
    	\int_{\mathbb{R}}V_{\lam}(\alpha)d\alpha =J_{k,s}(\bm{\lambda}) P^{s-k}.
    \end{align} 
 We complete the proof of Lemma \ref{major arc} by combining \eqref{err4} and \eqref{err5}.  \end{proof}

    \section{The minor arc}\label{sec:minor}
    The task in this section is to prove Lemma \ref{minor arc}.
  We first establish an auxiliary lemma.
    \begin{lemma}\label{lem1}
    	Let $h,\Delta>0$ and $s>1$. Suppose that 
    	\begin{align*}
    		\Theta=\left\lbrace\theta_1,\theta_2,\cdots,\theta_m \right\rbrace \subset \mathbb{R}
    	\end{align*}
    	satisfies
    	\begin{align*}
    		\left|\theta_i-\theta_j\right| >\Delta
    	\end{align*}
    	for any $i\neq j$. Then we have
    	\begin{align*}
    		\int_{1}^{2}\left(\sum_{i=1}^{m}\min \big(1,\,\frac{h}{\left|u-\theta_i\right|}\big)\right)^s du\ll h^s\Delta^{-s}\log^s (2m)+mh,
    	\end{align*}where the implicit constant depends only on $s$.
    \end{lemma}
    \begin{proof}Choose $1\leqslant i_0\leqslant m$ (depending on $u$) such that 
    $$|u-\theta_{i_0}|=\min_{1\leqslant i\leqslant m}|u-\theta_i|.$$
    Then one has\begin{align*}
    		\sum_{\substack{i=1 \\ i\not=i_0}}^{m}\min \big(1,\,\frac{h}{\left|u-\theta_i\right|}\big)\ll \sum_{j=1}^{m}\frac{h}{\Delta j}\ll h\Delta^{-1}\log(2m),
    	\end{align*}
   and thus, \begin{align*}
    		\left(\sum_{i=1}^{m}\min\big(1,\,\frac{h}{\left|u-\theta_i\right|}\big)\right)^s&\,\ll 
    \min\big(1,\,\frac{h}{\left|u-\theta_{i_0}\right|}\big)^s+h^s\Delta^{-s}\log^s (2m)
    \\ &\,\ll 
    \sum_{i=1}^{m}\min\big(1,\,\frac{h}{\left|u-\theta_{i}\right|}\big)^s+h^s\Delta^{-s}\log^s (2m)
    	\end{align*}
    Note that $\int_1^2\min\big(1,\,\frac{h}{\left|u-\theta\right|}\big)^s du\ll h$ for any $\theta\in \mathbb{R}$. Now  we deduce that
       	\begin{align*}
    		\int_{1}^{2}\left(\sum_{i=1}^{m}\min \big(1,\frac{h}{\left|u-\theta_i\right|}\big)\right)^s du&\,\ll \sum_{i=1}^{m}\int_1^2\min\big(1,\,\frac{h}{\left|u-\theta_{i}\right|}\big)^s du+h^s\Delta^{-s}\log^s (2m)
    \\ &\,\ll mh+h^s\Delta^{-s}\log^s (2m).
    	\end{align*}
This completes the proof of this lemma.
    \end{proof}
    
We introduce 
      	\begin{align*}
    		\mathcal{J}(P)=\int_{Q}^{+\infty}\int_{\alpha}^{+\infty}\left(\sum_{x\leqslant P}\sum_{y\leqslant P}\min\big(1,\,\frac{1}{\left|\alpha x^k-\beta y^k\right|}\big)\right)^sK_{\eta}(\alpha)K_{\eta}(\beta) d\beta d\alpha,
    	\end{align*}
    where the notation $x\leqslant P$ means $1\leqslant x\leqslant P$. One has the following.   
\begin{lemma}
\label{lemmaJP}
    	Suppose that $0< \eta\leqslant 2$. We have
    	\begin{align*}
    		\mathcal{J}(P)\ll \Big(\eta^2P^{2s-2k-\frac{2s-2k}{k}}+\eta P^{s-k}\Big)\log^{4s}P.
    	\end{align*}
    \end{lemma}    
 \begin{proof}We change variables by $\beta=\alpha 2^{rk}t^k$ with $r=0,1,2,\cdots$ and $1\leqslant t<2$ to deduce that
        \begin{align*}
    		\mathcal{J}(P)=\int_{Q}^{+\infty}\alpha K_{\eta}(\alpha)\sum_{r=0}^{\infty}2^{rk}\int_{1}^{2}kt^{k-1}K_{\eta}(\alpha2^{rk}t^{k})
    \left(\sum_{x\leqslant P}\sum_{y\leqslant P}\min\big(1,\,\frac{1}{\alpha\left| x^k-2^{rk}t^k y^k\right|}\big)\right)^sdtd\alpha.
        \end{align*}
        
        We introduce
        \begin{align*}
        	\Phi(U,V):=\Phi_{\alpha,r,t}(U,V)=\sum_{U\leqslant x<2U}\sum_{V\leqslant y<2V}\min\big(1,\,\frac{1}{\alpha\left| x^k-2^{rk}t^k y^k\right|}\big)
        \end{align*}
By the dyadic argument, 
        \begin{align*}\sum_{x\leqslant P}\sum_{y\leqslant P}\min\big(1,\,\frac{1}{\alpha\left| x^k-2^{rk}t^k y^k\right|}\big)\leqslant \sum_{U,V}\Phi(U,V),
        \end{align*}
        where $U,V$ are taken over powers of $2$ no more than $P$.
        Then by H\"older's inequality,
                \begin{align}\label{JtoL}
    		\mathcal{J}(P)\ll (\log P)^{2s}\sup_{U,V\leqslant P}\mathcal{L}(U,V),\end{align}
    where  \begin{align}\label{defineL}
    \mathcal{L}(U,V)=\int_{Q}^{\infty}\alpha K_{\eta}(\alpha)\sum_{r=0}^{\infty}2^{rk}\int_{1}^{2}kt^{k-1}K_{\eta}(\alpha2^{rk}t^{k})
    \Phi(U,V)^sdtd\alpha.
        \end{align}

      \noindent {\bf Case 1.} $\frac{1}{8}U>2^{r}V\;\text{or}\;8U<2^{r}V$.  In this case, we have 
        \begin{align*}
        	\left|x^k-2^{rk}t^ky^k\right|\gg \max\big(U,2^{r}V\big)^k\geqslant (2^{r}UV)^{\frac{k}{2}}.
        \end{align*}
        Then we deduce that\begin{align*}
    \Phi(U,V)\ll UV\min\big(1,\frac{1}{\alpha(2^{r}UV)^{\frac{k}{2}}}\big)=&\,\min\big( UV,\frac{ UV}{\alpha(2^{r}UV)^{\frac{k}{2}}}\big)
    \\ \ll &\,  (UV)^{1-\frac2k}\Big(\frac{ UV}{\alpha(2^{r}UV)^{\frac{k}{2}}}\Big)^{\frac{2}{k}}=\alpha^{-\frac{2}{k}}2^{-r}.
        \end{align*}
 Inserting this into \eqref{defineL}, we have  
        \begin{align*} \mathcal{L}(U,V)&\,\ll\int_{Q}^{\infty}\alpha K_{\eta}(\alpha)\sum_{r=0}^{\infty}2^{rk}\int_{1}^{2}K_{\eta}(\alpha2^{rk}t^{k})\alpha^{-\frac{2s}{k}}2^{-rs}dtd\alpha.
        \end{align*}
        Using $K_{\eta}(\cdot)\ll \eta$, we further obtain
       \begin{align}\label{boundLUV1} \mathcal{L}(U,V)\ll \eta^2\int_{Q}^{\infty}\alpha^{1-\frac{2s}{k}}d\alpha
       \ll \eta^2P^{2s-2k-\frac{2s-2k}{k}}.
        \end{align}  
 
         \noindent {\bf Case 2.}  $\frac{1}{8}U\leqslant 2^{r}V\leqslant 8U$. In this case, we have 
        \begin{align*}
        	|x^k-2^{rk}t^ky^k|\gg |\frac{x}{2^ry}-t|2^{rk}V^k.
        \end{align*} 
        Then       \begin{align*}
        	\Phi(U,V)\ll \sum_{U\leqslant x<2U}\sum_{V\leqslant y<2V}\min\big(1,\frac{1}{\alpha|\frac{x}{2^ry}-t|2^{rk}V^k}\big)
        \end{align*}
        We write $d=(x,y)$ and $x=dA,y=dB$. Then $(A,B)=1$. By splitting the sum over $d$ into dyadic intervals, we deduce that
 %       \begin{align*}
%        	\Phi(U,V)&\,\ll \sum_{d\leqslant 2^{1-r}P}d \sum_{\substack{A<2^{4+r}V/d\\ B<2V/d\\(A,B)=1}}\min(1,\frac{1}{\alpha\left|\frac{A}{2^rB}-t|2^{rk}V^k\right| })
%         \\ &\,\ll \sum_{D=2^j\leqslant 2^{1-r}P}D\Psi(D),
%        \end{align*}
       \begin{align*}
        	\Phi(U,V)\ll \sum_{d\leqslant 2V} \sum_{\substack{A<2^{4+r}V/d\\ B<2V/d\\(A,B)=1}}\min\big(1,\frac{1}{\alpha|\frac{A}{2^rB}-t|2^{rk}V^k }\big)
      \ll \sum_{D=2^j\leqslant 2V}D\Psi_t(D),
        \end{align*}
        where      \begin{align*}
        	\Psi_t(D)= \sum_{\substack{A<2^{5+r}V/D\\ B<4V/D\\(A,B)=1}}\min\big(1,\frac{1}{\alpha|\frac{A}{2^rB}-t|2^{rk}V^k }\big).
        \end{align*}
      Inserting above into \eqref{defineL} and applying H\"older's inequality, we obtain
      \begin{align*}
    \mathcal{L}(U,V)\ll (\log P)^{s}\sup_{D\leqslant 2V}D^{s}\int_{Q}^{\infty}\alpha K_{\eta}(\alpha)\sum_{r=0}^{\infty}2^{rk}\int_{1}^{2}kt^{k-1}K_{\eta}(\alpha2^{rk}t^{k})
    \Psi_t(D)^sdtd\alpha.
        \end{align*}

        By \eqref{propK1}, we have
        \begin{align*}
        	K_{\eta}(\alpha2^{rk}t^k)\ll \frac{\eta}{1+\left(\eta\alpha2^{rk}\right)^2}.
        \end{align*}
Then we obtain       \begin{align*}\int_{1}^{2}kt^{k-1}K_{\eta}(\alpha2^{rk}t^{k})
    \Psi_t(D)^sdt \ll \frac{\eta}{1+\left(\eta\alpha2^{rk}\right)^2}\int_{1}^{2}
    \Psi_t(D)^sdt.
        \end{align*} Therefore,
              \begin{align}\label{boundLUV2}
    \mathcal{L}(U,V)\ll (\log P)^{s}\sup_{D\leqslant 2V}D^{s}\int_{Q}^{\infty}\alpha K_{\eta}(\alpha)\sum_{r=0}^{\infty}2^{rk}\frac{\eta}{1+\left(\eta\alpha2^{rk}\right)^2}\int_{1}^{2}
    \Psi_t(D)^sdtd\alpha.
        \end{align} 
        
             The numbers $\frac{A}{2^rB}$ subject to $A<2^{5+r}V/D, B<4V/D, (A,B)=1$ are pairwise distinct and 
             $$\#\{\frac{A}{2^rB}:\ A<2^{5+r}V/D, B<4V/D, (A,B)=1\}\ll 2^rV^2D^{-2}.$$
   Moreover, their mutual space is
        $\geqslant \frac{D^2}{2^{r+4}V^2}$. By Lemma \ref{lem1} with $m=O(2^rV^2D^{-2})$, $\Delta=\frac{D^2}{2^{r+4}V^2}$ and $h=\frac{1}{\alpha 2^{rk}V^k}$,  we obtain
        \begin{align*}
        	\int_{1}^{2}
    \Psi_t(D)^sdt \ll 2^{rs(1-k)}V^{s(2-k)}D^{-2s}\alpha^{-s}(\log P)^s+2^{r(1-k)}V^{2-k}D^{-2}\alpha^{-1}.
        \end{align*}
 One also has the trivial bound 
        \begin{align*}
        	\int_{1}^{2}
    \Psi_t(D)^sdt \ll \frac{2^{rs}V^{2s}}{D^{2s}}.
        \end{align*}
        Therefore, 
        \begin{align*}
        	\int_{1}^{2}
    \Psi_t(D)^sdt \ll \min\Big(\frac{2^{rs}V^{2s}}{D^{2s}},\,2^{rs(1-k)}V^{s(2-k)}D^{-2s}\alpha^{-s}(\log P)^s+2^{r(1-k)}V^{2-k}D^{-2}\alpha^{-1}\Big).
        \end{align*}
        Multiplying by $D^s$ and noting that $D\leqslant 2V\ll 2^{-r}P$, we have
        \begin{align*}
        	D^s\int_{1}^{2}
    \Psi_t(D)^sdt \ll &\, \min\Big(2^{rs}V^{2s}D^{-s},\, 2^{rs(1-k)}V^{s(2-k)}D^{-s}\alpha^{-s}(\log P)^s+2^{r(1-k)}V^{s-k}\alpha^{-1}\Big)
    \\ \ll &\, \min\Big(2^{rs}V^{2s},\, 2^{rs(1-k)}V^{s(2-k)}\alpha^{-s}(\log P)^s+2^{r(1-s)}P^{s-k}\alpha^{-1}\Big).
        \end{align*}
        By using the inequality $\min (a,b+c)\leqslant a^{1-2/k}b^{2/k}+c$ for $a,b,c>0$, we further obtain that 
   \begin{align*}
        	D^s\int_{1}^{2}
    \Psi_t(D)^sdt \ll 2^{-rs}\alpha^{-\frac{2s}{k}}(\log P)^{2s/k}+2^{r(1-s)}P^{s-k}\alpha^{-1}.
        \end{align*} 
 Putting the above into \eqref{boundLUV2}, we conclude that
        \begin{align}\label{boundLUV3}
          \mathcal{L}(U,V)\ll (\log P)^{2s}\left(\int_{Q}^{+\infty}\eta K_{\eta}(\alpha)\alpha^{-\frac{2s}{k}+1}d\alpha
          +\int_{0}^{+\infty}K_{\eta}(\alpha)\sum_{r=0}^{+\infty}\frac{2^{r(1+k-s)}P^{s-k}\eta}{1+(\eta\alpha2^{rk})^2}d\alpha\right),
        \end{align} 
where we have used the condition $s>k$. Now the first integration in \eqref{boundLUV3} is
      \begin{align*}
      	\ll \eta^2P^{2s-2k-\frac{2s-2k}{k}}.
      \end{align*}  Interchanging the order of integration and summation, and then making the substitution $u=\eta\alpha 2^{rk}$, we find that the second integration in \eqref{boundLUV3} is
      \begin{align*}
      	P^{s-k}\sum_{r=0}^{+\infty}2^{r(1-s)}\int_{0}^{+\infty}K_{\eta }(\frac{u}{2^{rk}\eta})\frac{1}{1+u^2}du\ll \eta P^{s-k}.
      \end{align*}
      Therefore, we conclude that
              \begin{align}\label{boundLUV4}
          \mathcal{L}(U,V)\ll \Big(\eta^2P^{2s-2k-\frac{2s-2k}{k}}+\eta P^{s-k}\Big)\log^{2s}P.
        \end{align} 
        
   In view of \eqref{boundLUV1} proved in Case 1 and \eqref{boundLUV4} obtained in Case 2, we see that indeed \eqref{boundLUV4} holds in both cases. Then combining \eqref{JtoL} and \eqref{boundLUV4},      
      we complete the proof of Lemma \ref{lemmaJP}.
    \end{proof}

 Now we are able to establish Lemma \ref{minor arc}. 
 
 \begin{proof}[Proof of Lemma \ref{minor arc}] Note that
    	\begin{align*}
    		\left|I(\lam;\mathfrak{m})\right|^2=\int_{\mathfrak{m}}\int_{\mathfrak{m}}f(\lambda_1\alpha)\overline{f(\lambda_1\beta)}\cdots f(\lambda_s\alpha)\overline{f(\lambda_s\beta)}K_{\eta}(\alpha)K_{\eta}(\beta)d\alpha d\beta.
    	\end{align*}
 We deduce that
    	\begin{align*}
    		\int_{\la(C)}\left|I(\lam;\mathfrak{m})\right|^2d\lam &\,\leqslant \int_{[-C,C]^s}\int_{\mathfrak{m}}\int_{\mathfrak{m}}f(\lambda_1\alpha)\overline{f(\lambda_1\beta)}\cdots f(\lambda_s\alpha)\overline{f(\lambda_s\beta)}K_{\eta}(\alpha)K_{\eta}(\beta)d\alpha d\beta d\lam.
 \\ &\,=\int_{\mathfrak{m}}\int_{\mathfrak{m}}\left(\int_{-C}^{C}f(\lambda\alpha)\overline{f(\lambda\beta)}d\lambda\right)^sK_{\eta}(\alpha)K_{\eta}(\beta) d\alpha d\beta.
    	\end{align*}
    	By the definition of $f$ in \eqref{definef}, we have
    	\begin{align*}
    		\int_{-C}^{C}f(\lambda\alpha)\overline{f(\lambda\beta)}d\lambda=\sum_{\left|x\right|\leqslant P}\sum_{\left|y\right|\leqslant P}\int_{-C}^{C}e\left(\lambda(\alpha x^k-\beta y^k)\right)d\lambda
    	\end{align*}
 Then we have
    	\begin{align*}
    		\int_{-C}^{C}f(\lambda\alpha)\overline{f(\lambda\beta)}d\lambda\ll \sum_{\left|x\right|\leqslant P}\sum_{\left|y\right|\leqslant P}\min\big(1,\frac{1}{\left|\alpha x^k-\beta y^k\right|}\big).
    	\end{align*}
    Thus, 
   \begin{align*}
    		\int_{\la(C)}\left|I(\lam;\mathfrak{m})\right|^2d\lam\ll \int_{\mathfrak{m}}\int_{\mathfrak{m}}\left(\sum_{\left|x\right|\leqslant P}\sum_{\left|y\right|\leqslant P}\min\big(1,\frac{1}{\left|\alpha x^k-\beta y^k\right|}\big)\right)^sK_{\eta}(\alpha)K_{\eta}(\beta) d\alpha d\beta.
    	\end{align*}

    	By noting that $|\alpha x^k-\beta y^k|\geqslant \left||\alpha||x|^k-|\beta||y|^k\right|$, we have 
      \begin{align*}
    		\int_{\la(C)}\left|I(\lam;\mathfrak{m})\right|^2d\lam\ll \int_{Q}^{+\infty}\int_{Q}^{+\infty}\left(\sum_{0\leqslant x\leqslant P}\sum_{0\leqslant y\leqslant P}\min\big(1,\frac{1}{\left|\alpha x^k-\beta y^k\right|}\big)\right)^sK_{\eta}(\alpha)K_{\eta}(\beta) d\alpha d\beta.
    	\end{align*}  
By the symmetry in $\alpha$ and $\beta$, the contribution from $x\geqslant 1$ and $y\geqslant 1$ is bounded by $\mathcal{J}(P)$. Thus, 
      \begin{align}\label{boundtoJPand}
    		\int_{\la(C)}\left|I(\lam;\mathfrak{m})\right|^2d\lam\ll \mathcal{J}(P)+\int_{Q}^{+\infty}\int_{Q}^{+\infty}\left(\sum_{0\leqslant y\leqslant P}
    \min\big(1,\frac{1}{\beta y^k}\big)\right)^sK_{\eta}(\alpha)K_{\eta}(\beta) d\alpha d\beta.
    	\end{align}
Since $\sum_{0\leqslant y\leqslant P}\min\big(1,\frac{1}{\beta y^k}\big)\ll 1+\beta^{-\frac{1}{k}}$, one has 
      \begin{align*}&\,\int_{Q}^{+\infty}\int_{Q}^{+\infty}\left(\sum_{0\leqslant y\leqslant P}\min\big(1,\frac{1}{\beta y^k}\big)\right)^sK_{\eta}(\alpha)K_{\eta}(\beta) d\alpha d\beta
      \\ \ll&\, \int_{Q}^{+\infty}\int_{Q}^{+\infty}K_{\eta}(\alpha)K_{\eta}(\beta) d\alpha d\beta+\int_{Q}^{+\infty}\int_{Q}^{+\infty}\beta^{-\frac{s}{k}}K_{\eta}(\alpha)K_{\eta}(\beta) d\alpha d\beta.
    	\end{align*}
By using \eqref{propK0} and \eqref{PropK2}, we further have    
    \begin{align}\label{and}\int_{Q}^{+\infty}\int_{Q}^{+\infty}\left(\sum_{0\leqslant y\leqslant P}\min\big(1,\frac{1}{\beta y^k}\big)\right)^sK_{\eta}(\alpha)K_{\eta}(\beta) d\alpha d\beta
\ll 1+\eta Q^{-\frac{s}{k}+1}\ll 1+\eta P^{s-k-\frac{s-k}{k}}.
    	\end{align}
    
    Note that $\eta P^{s-k-\frac{s-k}{k}}$ is bounded by $1+\eta^2P^{2s-2k-\frac{2s-2k}{k}}(\log P)^{4s}$. Now Lemma \ref{minor arc} follows from Lemma \ref{lemmaJP}, \eqref{boundtoJPand} and \eqref{and} immediately.\end{proof}

 \section{Proof of Theorem~\ref{thm2}}\label{sec:proof}

   	Put
 	\begin{align*}
 		E(P,\eta,\lam)=I(P,\eta,\lam)-\eta J_{k,s}(\bm{\lambda}) P^{s-k}.
 	\end{align*}
  	\begin{lemma}\label{7}
 			Let $s,k$ be integers with $s>k\geqslant 2$. Suppose that $\delta= \frac{1}{20k}$. For almost all $\lam\in \la $, there exists a number $T_0=T_0(\lam,k)$ such that the inequality 
 		\begin{align*}
 			\left|E(P,\eta,\lam)\right|\leqslant 3\eta P^{s-k-3\delta}
 		\end{align*}
 		holds for all $P\geqslant T_0$ and $\frac{1}{2}P^{-\frac{1}{2}}<\eta\leqslant 2$.
 	\end{lemma}
\begin{proof} 
For $z\in \{1,2\}$, we define
 	\begin{align*}
 		\mathscr{K}(P,\eta,z)=\left\lbrace \lam\in \la(C):\,\left|E(P,\eta,\lam)\right|>z\eta P^{s-k-3\delta}\right\rbrace ,
 	\end{align*}
 	and
 	\begin{align*}
 		\mathscr{K}(P,z)=\bigcup_{\frac{1}{4}P^{-\frac12}\leqslant \eta\leqslant 2}\mathscr{K}(P,\eta,z).
 	\end{align*}
 	Let $l_0$ be the smallest integer satisfying
 	\begin{align*}
 		\frac{1}{4}P^{-\frac12}\left(1+P^{-4\delta}\right)^{l_0}\geqslant 2.
 	\end{align*}
 	We define
 	\begin{align*}
 		\eta_{l}=\min \left(2,\ \frac{1}{4}P^{-\frac12}\big(1+P^{-4\delta}\big)^{l}\right)
 	\end{align*}
 	for $0\leqslant l\leqslant l_0$.
 	We have
 	\begin{align*}
 		l_0\ll P^{4\delta}\log P.
 	\end{align*}
 For $\frac{1}{4}P^{-\theta}\leqslant \eta\leqslant 2$, there exists $l~(0\leqslant l\leqslant l_0)$ such that $\eta_{l}\leqslant \eta< \eta_{l+1}$. Note that 
  $$\left|\eta_{l+1}-\eta_{l}\right|\ll \eta_{l} P^{-4\delta}.$$
For sufficiently large $P$, we claim that
\begin{align}\label{claim}\mathscr{K}(P,\eta,2)\subseteq \mathscr{K}(P,\eta_{l},1)\cup \mathscr{K}(P,\eta_{l+1},1).\end{align}
Otherwise, there exists $\lam$ such that 
$|I(P,\eta,\lam)-\eta J_{k,s}(\bm{\lambda}) P^{s-k}|>2\eta P^{s-k-3\delta}$ but 
\begin{align*}|I(P,\eta_l,\lam)-\eta_l J_{k,s}(\bm{\lambda}) P^{s-k}|\leqslant\eta_l P^{s-k-3\delta} \ \textrm{ and }
\ \ |I(P,\eta_{l+1},\lam)-\eta_{l+1}J_{k,s}(\bm{\lambda}) P^{s-k}|\leqslant\eta_{l+1} P^{s-k-3\delta}.\end{align*}
We deduce from above that
\begin{align*}I(P,\eta,\lam)-\eta J_{k,s}(\bm{\lambda}) P^{s-k}&\,\leqslant I(P,\eta_{l+1},\lam)-\eta_{l+1}J_{k,s}(\bm{\lambda}) P^{s-k}+(\eta_{l+1}-\eta)J_{k,s}(\bm{\lambda}) P^{s-k}
\\ &\, \leqslant \eta  P^{s-k-3\delta}+O(\eta J_{k,s}(\bm{\lambda}) P^{s-k-4\delta}).\end{align*}
Similarly, 
\begin{align*}I(P,\eta,\lam)-\eta J_{k,s}(\bm{\lambda}) P^{s-k}\geqslant -\eta  P^{s-k-3\delta}+O(\eta J_{k,s}(\bm{\lambda}) P^{s-k-4\delta}).\end{align*}
Now for large $P$, one has $|I(P,\eta,\lam)-\eta J_{k,s}(\bm{\lambda}) P^{s-k}|<2\eta P^{s-k-3\delta}$, which is a contradiction. Thus, the claim \eqref{claim} is true.

 	We conclude from \eqref{claim} that
 	\begin{align*}
 		\mathscr{K}(P,2)\subseteq\bigcup_{l=0}^{l_0}\mathscr{K}(P,\eta_{l},1).
 	\end{align*}
 	Hence we have
 	\begin{align}\label{boundKP2}
 		\int_{\mathscr{K}(P,2)}d\lam\leqslant \sum_{l=0}^{l_0}\int_{\mathscr{K}(P,\eta_{l},1)}d\lam.
 	\end{align}
Lemmas~\ref{major arc}
 	and~\ref{minor arc} give
 	\begin{align*}
 		\int_{\la(C)}|E(P,\eta,\lam)|^{2}\,d\lam
 		\ll_{C}(\log P)^{4s}
 		\left(1+\eta P^{s-k}
 		+\eta^{2}P^{2s-2k-\frac{1}{k}}\right)
 	\end{align*}
    for any $0<\eta\leqslant  2$. Thus, we have
%    \begin{align*}
%    	\int_{\mathscr{K}(P,\eta_{l},1)}d\lam\ll &\, \Big(\eta_{l}^{-2}P^{6\delta-2(s-k)}+\eta_{l}^{-1}P^{6\delta-(s-k)}+P^{6\delta-\frac1k}\Big)(\log P)^{3s}
%  \\ \ll &\,  P^{6\delta-\frac{1}{k}}(\log P)^{3s}.
%    \end{align*}
    \begin{align*}
    	\int_{\mathscr{K}(P,\eta_{l},1)}d\lam\ll  \Big(\eta_{l}^{-2}P^{6\delta-2(s-k)}+\eta_{l}^{-1}P^{6\delta-(s-k)}+P^{6\delta-\frac1k}\Big)(\log P)^{4s}
 \ll \,  P^{6\delta-\frac{1}{k}}(\log P)^{4s}.
    \end{align*}
 We put the above into \eqref{boundKP2} and obtain
 	\begin{align}\label{measure}
 		\int_{\mathscr{K}(P,2)}d\lam\ll  P^{10\delta-\frac{1}{k}}(\log P)^{4s+1}
 	\end{align}
 for all sufficiently large $P$. 
 
 Let $P_j=j^{\frac{1}{4\delta}}$. Suppose that $J$ is a large integer such that \eqref{measure} holds for all $P\geqslant P_J$. Thus, by \eqref{measure}, we have
		\begin{align*}
			m\left(\bigcup_{j=J}^{\infty}\mathscr{K}(P_j,2)\right)\leqslant \sum_{j=J}^{\infty}\int_{\mathscr{K}(P_j,2)}d\lam <\infty.
		\end{align*}
		We write
		\begin{align*}
			\mathscr{K}=\bigcap_{i=J}^{\infty}\bigcup_{j=i}^{\infty}\mathscr{K}(P_j,2).
		\end{align*}
		By the Borel-Cantelli Lemma, the Lebesgue measure of $\mathscr{K}$ is zero. For any $\lam \in \la (C)\setminus\mathscr{K}$, there exists a positive integer $J_{\lam}$, such that
		\begin{align}\label{aaa}
			\lam \notin \bigcup_{j=J_{\lam}}^{\infty}\mathscr{K}(P_j,2).
		\end{align}
	For $P\geqslant P_{J_{\lam}}$, there is a positive integer $j\geqslant J_{\lam}$ such that $P_j\leqslant P<P_{j+1}$, and then
		\begin{align*}
			\max\left(\left|P-P_j\right|,\left|P-P_{j+1}\right|\right)\ll\left|P_j-P_{j+1}\right|\ll P_j^{1-4\delta} \ll P^{1-4\delta}.
		\end{align*}
One has
\begin{align}\label{P}
			P_{j}^{s-k}-P^{s-k},\, P_{j+1}^{s-k}-P^{s-k}\ll P^{s-k-4\delta}.
		\end{align}
The function $I(P,\eta,\lam)$ is non-decreasing as a function of $P$. Thus, we have
		\begin{align}\label{I}
			I(P_j,\eta,\lam )\leqslant I(P,\eta,\lam)\leqslant I(P_{j+1},\eta,\lam).
		\end{align}
		By the definition of $\mathscr{K}(P,2)$, we have
		\begin{align}\label{II}
			\begin{aligned}
			I(P_j,\eta,\lam)\geqslant \eta \J P_j^{s-k}-2\eta P_j^{s-k-3\delta},\\I(P_{j+1},\eta,\lam)\leqslant \eta \J P_{j+1}^{s-k}+2\eta P_{j+1}^{s-k-3\delta}
		    \end{aligned}
		\end{align}
	 for $\frac{1}{3}P^{-\frac{1}{2}}\leqslant \eta\leqslant 2$.
		Combining \eqref{P}, \eqref{I} and \eqref{II}, we obtain
		\begin{align*}
			|E(P,\eta,\lam)|\leqslant 2\eta P^{s-k-3\delta}+O(\eta P^{s-k-4\delta}).
		\end{align*}Since $P$ is sufficiently large, we further have
		\begin{align*}
			|E(P,\eta,\lam )|\leqslant 3\eta P^{s-k-3\delta}
		\end{align*}
		whenever $\frac{1}{2}P^{-\frac{1}{2}}\leqslant \eta\leqslant 2$.

	We call $\lam$ exceptional if it does not satisfy the statement of Lemma~\ref{7}. Then, the above argument gives that the exceptional $\lam\in\la(C)$ must belong to $\mathscr{K}$. Thus the exceptional set also has Lebesgue measure zero. We consider the union of these sets, with $C\geqslant 2$ running over natural numbers, to complete the proof of Lemma~\ref{7}.
\end{proof}

 	We now turn to removing the weights from the counting function $I(P,\eta,\lam)$ by a sandwich technique in \cite[Section 8]{BD}.

	\begin{theorem}\label{thm3}
		Let $s,k$ be integers with $s>k\geqslant 2$. Suppose that $\delta= \frac{1}{20k}$.
		Then for almost all $\bm{\lambda}\in\la$,
		there exists a number
		$T_0=T_0(\bm{\lambda},k)$
		and a positive real number
		$J_{k,s}(\bm{\lambda})$
		such that the inequality
		\begin{align*}
			\left|
			N_{\bm{\lambda}}^{(k)}(P;\tau)
			-2\tau J_{k,s}(\bm{\lambda}) P^{s-k}
			\right|
			<\tau P^{s-k-\delta}
		\end{align*}
		holds for all $P\geqslant T_0$ and 
		 $P^{-\frac{1}{2}}\leqslant\tau\leqslant 1$.
	\end{theorem}
\begin{proof}
	Let $P^{-\frac{1}{2}}\leqslant \tau\leqslant 1$, and put $\Delta=\tau P^{-\frac{3}{2}\delta}$. Following \cite[Section 8]{BD}, we introduce
	\begin{align}\label{W}
		W^{\pm}(\alpha)=\pm \Delta^{-1}\left(\left(\tau\pm \Delta\right)W_{\tau\pm \Delta}(\alpha)-\tau W_{\tau}(\alpha)\right).
	\end{align}
Note that $W^{\pm}$ are continuous functions with $0\leqslant W^{\pm}(\alpha)\leqslant 1$ for all $\alpha\in\mathbb{R}$ and 
	\begin{align*}
		W^{+}(\alpha)=1\;\;\text{for}\;\left|\alpha\right|\leqslant \tau,\;\;\;\;\;\;\;\;\;\;\;\;&W^{+}(\alpha)=0\;\;\text{for}\;\left|\alpha\right|\geqslant \tau+\Delta,\\W^{-}(\alpha)=1\;\;\text{for}\;\left|\alpha\right|\leqslant \tau-\Delta,\;\;\;\;\;&W^{-}(\alpha)=0\;\;\text{for}\;\left|\alpha\right|\geqslant \tau.
	\end{align*}
    Now define the counting functions
    \begin{align*}
    	I^{\pm}(P,\tau,\lam)=\sum_{\left|\x\right|_{\infty}\leqslant P}W^{\pm}\left(F_{\lam}(\x)\right),
    \end{align*}
    and one has
    \begin{align*}\label{sand}
    	I^{-}(P,\tau,\lam)\leqslant N_{\bm{\lambda}}^{(k)}(P;\tau)\leqslant I^{+}(P,\tau,\lam).
    \end{align*}
    By \eqref{W}, we have
    \begin{align*}
    	I^{\pm}(P,\tau,\lam)=\pm \Delta^{-1}\left(\left(\tau\pm \Delta\right)I(P,\tau\pm \Delta,\lam )-\tau I(P,\tau,\lam)\right).
    \end{align*}
    Now suppose that $\lam$ is not in the exceptional null set excluded from the claim in Lemma~\ref{7}. Then one has
    \begin{align*}
    	I^{\pm}(P,\tau,\lam)=&\, \pm \Delta^{-1}\left(\big((\tau\pm \Delta)^2-\tau^2\big)J_{k,s}(\bm{\lambda}) P^{s-k}+O(\tau P^{s-k-\frac{3}{2}\delta})\right)
    \\ =&\, 2\tau J_{k,s}(\bm{\lambda}) P^{s-k}+O(\tau P^{s-k-\frac{3}{2}\delta}).
    \end{align*}
    We conclude from above that
      \begin{align*}
    	N_{\bm{\lambda}}^{(k)}(P;\tau)= 2\tau J_{k,s}(\bm{\lambda}) P^{s-k}+O(\tau P^{s-k-\frac{3}{2}\delta}).
    \end{align*} 
    Then for $P$ sufficiently large, we obtain the inequality 
        \begin{align*}
    	\left|N_{\bm{\lambda}}^{(k)}(P;\tau)-2\tau J_{k,s}(\bm{\lambda}) P^{s-k}\right|<\tau P^{s-k-\delta}
    \end{align*}
    for almost all $\lam$. This completes the proof of Theorem \ref{thm3}.
\end{proof}

\begin{proof}[Proof of Theorem~\ref{thm2}] In view of Theorem \ref{thm3}, it suffices to consider $0<\tau<P^{-\frac{1}{2}}$. Suppose that $\lam$ is not in the exceptional set claimed in Theorem \ref{thm3}. Then for $P\geqslant T_0$ with $T_0$ in Theorem \ref{thm3}, one has
$$N_{\bm{\lambda}}^{(k)}(P;\tau)\leqslant N_{\bm{\lambda}}^{(k)}(P;P^{-\frac{1}{2}})\leqslant 2P^{-\frac{1}{2}}J_{k,s}(\bm{\lambda}) P^{s-k}+P^{s-k-\delta-\frac{1}{2}}\ll P^{s-k-\frac{1}{2}}.$$
Thus, 
$$\left|N_{\bm{\lambda}}^{(k)}(P;\tau)-2\tau J_{k,s}(\bm{\lambda}) P^{s-k}\right|\leqslant N_{\bm{\lambda}}^{(k)}(P;\tau)+2\tau J_{k,s}(\bm{\lambda}) P^{s-k}\ll P^{s-k-\frac{1}{2}}.$$
Then there exists $T_0'$ such that $\left|N_{\bm{\lambda}}^{(k)}(P;\tau)-2\tau J_{k,s}(\bm{\lambda}) P^{s-k}\right|< P^{s-k-\delta}$ for all $P\geqslant T_0'$. 
We complete the proof of Theorem \ref{thm2} by taking $P_0=\max(T_0,\,T_0')$.
\end{proof}
\noindent {\bf Acknowledgement}.
	This work is supported by the National Key Research and Development Program of China (No. 2021YFA1000701) and National Natural Science Foundation of China (No. 12471088).

\end{document}